\documentclass{article}
\usepackage[latin1]{inputenc}
\usepackage{amsmath,euscript}
\usepackage{amsthm}
\usepackage{amsfonts}

\theoremstyle{plain}
\newtheorem{lemma}{Lemma}[section]
\newtheorem{proposition}[lemma]{Proposition}
\newtheorem{theorem}[lemma]{Theorem}
\newtheorem{corollary}[lemma]{Corollary}

\theoremstyle{remark}
\newtheorem*{remark}{Remark}
\newtheorem*{example}{Example}

\renewcommand{\phi}{\varphi}
\renewcommand{\epsilon}{\varepsilon}
\renewcommand{\rho}{\varrho}
\renewcommand{\d}{\boldsymbol{\delta}}
\renewcommand{\r}{\boldsymbol{\rho}}
\newcommand{\sprod}[2]{\langle#1,#2\rangle}

\begin{document}

\title{On the Non-Equivalence of Rearranged Walsh and Trigonometric
  Systems in $L_p$}
\author{Aicke Hinrichs \ and
  J{\"o}rg Wenzel\thanks{The authors were supported by DFG grants Hi
    584/2-2 and We 1868/1-1, resp.
    \newline Both authors thank the University of Pretoria where part
    of this work was carried out.
    }}

\maketitle

\begin{center}
  {\it Dedicated to Professor A. Pe{\l}czy{\'n}ski on the occasion of
    his $70^{th}$ birthday}
\end{center}

\begin{abstract}
  We consider the question whether the trigonometric system can be
  equivalent to some rearrangement of the Walsh system in $L_p$ for
  some $p\neq 2$. We show that this question is closely related to a
  combinatorial problem. This enables us to prove non-equivalence for
  a number of rearrangements. Previously this was known for the
  Walsh-Paley order only.
\end{abstract}

%%%%%%%%%%%%%%%%%%%%%%%%%%%%%%%%%%%%%%%%%%%%%%%%%%%%%%%%%%%%%%%%%%%%%%%

\section{Introduction}
\label{sec:intro}

Both the Walsh system and the trigonometric system are systems of
characters on a compact abelian group. This explains that many of the
results in the theory of those systems are parallel. However, those
similarities do usually not extend to the case when the systems are
compared directly. So it is known that the Walsh system in the
Walsh-Paley order and the trigonometric system are not equivalent in
$L_p$ for $p\neq 2$, see \cite{young76:_walsh_fourier}. A
``power-type" non-equivalence for those systems was recently shown in
\cite{wojtaszczyk00:_non_walsh}.

It does not seem natural to fix the order of the systems in this basis
equivalence problem. In \cite{wojtaszczyk00:_non_walsh} the conjecture
was made that non-equivalence also holds for arbitrary rearrangements
of the Walsh system. Nevertheless, the methods used in that paper are
very particular to the case of the Walsh-Paley order. The aim of this
note is to address the more general equivalence problem.

In a first part, we relate the equivalence question for a fixed
ordering to a question of algebraic combinatorial type. In a second
part, we apply this approach to prove non-equivalence for a number of
orderings. We obtain estimates of power type but we do not attempt to
find the optimal estimates here.

The method can be generalized to deal with the equivalence of
arbitrary systems of characters on compact abelian groups. This will
be studied elsewhere.

We consider the trigonometric system $(e_n)_{n\in\mathbb{Z}}$ on
$[0,1]$ given by the functions $e_n(t)= \exp(2 \pi i n t)$. Let
$r_0,r_1,\ldots$ be the system of Rademacher functions on $[0,1]$. For
$n=0,1,\ldots$ the binary expansion of $n$ is $n=\sum_{i=0}^\infty n_i
2^i$ with $n_i\in\{0,1\}$. Observe that the sum is a finite sum. The
$n$-th Walsh function is then given by
\[ w_n = \prod_{i=0}^\infty r_i^{n_i}.
\]

Fix a permutation $\sigma$ of $\{0,1,\ldots\}$. Given $1\le p\le
\infty$, we say that the trigonometric system and the Walsh system
rearranged with $\sigma$ are equivalent in $L_p$ if there exists $c>0$
such that the inequalities
\[ \frac{1}{c} \Big\| \sum_{k=0}^\infty \xi_k e_k \Big\|_p \le
\Big\| \sum_{k=0}^\infty \xi_k w_{\sigma(k)} \Big\|_p \le
c \Big\| \sum_{k=0}^\infty \xi_k e_k \Big\|_p
\]
hold for all sequences $(\xi_k)$ of complex numbers with only finitely
many nonzero terms. Here $\| \cdot \|_p$ denotes the norm in
$L_p[0,1]$.  We then write $(e_k)_{k\ge 0} \sim_p
(w_{\sigma(k)})_{k\ge 0}$. Observe that we actually only consider one
half of the trigonometric system. This is not really essential for
what follows but it simplifies the exposition significantly.

We now describe the organization of the paper in some detail. The next
section provides a basic duality result which shows that the
equivalence questions in $L_p$ and $L_{p'}$ are essentially the same
where as usual $p'$ denotes the conjugate number of $p$ given by
$1/p+1/p'=1$.

In Section 3 we introduce and study a sequence of functions crucial
for our purpose. Norm estimates for these functions provide a tool to
prove non-equivalence of Walsh and trigonometric systems. This is a
generalization of the method used in \cite{wojtaszczyk00:_non_walsh}.
Moreover, it will turn out to be fundamental for our considerations
that non-equivalence of Walsh and trigonometric systems in $L_p$ for
all $p\neq 2$ can be derived from a nontrivial $L_{p_0}$-norm estimate
of those functions for one fixed $p_0>2$.

In Section 4 we show that the $L_4$-norm of the key functions is
determined by the solution of a particular combinatorial problem. This
is due to the fact that the fourth power of the $L_4$-norm is a
polynomial function and to the orthogonality of Walsh and
trigonometric functions.

Section 5 applies this approach to concrete rearrangements of the
Walsh system, in particular to linear and piecewise linear
rearrangements and ``small" perturbations thereof. This includes all
the commonly used orderings of the Walsh functions, i.e. the
Walsh-Paley order, the original Walsh order, the Walsh-Kaczmarz order,
and the Walsh-Kronecker orders. We give the definitions of these
orderings at the appropriate places. More information on Walsh
functions and Walsh series can be found in the monographs
\cite{golubov91:_walsh_series} and \cite{schipp90:_walsh}.

%%%%%%%%%%%%%%%%%%%%%%%%%%%%%%%%%%%%%%%%%%%%%%%%%%%%%%%%%%%%%%%%%%%%%%%

\section{Duality}

To study the equivalence problem we introduce some notation. By
$\EuScript A_n$ we denote a general orthonormal system
$(a_0,\dots,a_{n-1})$ in $L_2[0,1]$, usually this will be the system
$\EuScript W_n$ of the first $n$ Walsh functions $(w_0,\dots,w_{n-1})$
or the system $\EuScript E_n$ of the first $n$ exponential functions
$(e_0,\dots,e_{n-1})$. Given a finite set
$\mathbb{F}\subseteq\{0,\dots,n-1\}$, we denote by $\EuScript
A(\mathbb{F})$ the system formed by the functions $a_k$ with
$k\in\mathbb{F}$. If in particular $\mathbb{F}=[m]:=\{0,\dots,m-1\}$
is the set of the first $m$ members of $\{0,\dots,n-1\}$ we write
again $\EuScript A_m$ for $\EuScript A([m])$. Given a permutation
$\sigma$ of $\{0,1,\dots\}$ we denote by $\EuScript
A^\sigma(\mathbb{F})$ the system formed by all $a_{\sigma(k)}$, where
$k\in \mathbb{F}$. Note that $\EuScript A^\sigma(\mathbb{F})$ and
$\EuScript A(\sigma(\mathbb{F}))$ differ just by their order.

We let $\r_p(\EuScript A(\mathbb{F}),\EuScript B(\mathbb{F}))$ denote
the smallest constant $c$ such that
\[ \Big\| \sum_{k\in\mathbb{F}} \xi_k a_k \Big\|_p \leq c
\Big\| \sum_{k\in\mathbb{F}} \xi_k b_k \Big\|_p
\]
holds for all complex numbers $\xi_k$ with $k\in\mathbb{F}$.

We are interested in the quantities
\begin{eqnarray*}
  \r_p(\EuScript E_n,\EuScript W^\sigma_n), & \hspace{2em} &
  \r_p(\EuScript W_n,\EuScript E^{\sigma^{-1}}_n), \\
  \r_p(\EuScript W^\sigma_n,\EuScript E_n), & \hspace{2em} &
  \r_p(\EuScript E^{\sigma^{-1}}_n,\EuScript W_n),
\end{eqnarray*}
where $\sigma^{-1}$ is the inverse permutation of $\sigma$.

In order to get some information on duality, we need another quantity,
which behaves better under passing from $p$ to $p'$. The $k$-th {\em
  Fourier coefficient} of a function $f \in L_p$ with respect to the
system $\EuScript B_n$ is given by $\sprod f{b_k} = \int_0^1 f(t)
\overline{b_k(t)} \, dt$. For simplicity, we will henceforth assume
that all orthonormal systems considered consist of bounded functions,
so Fourier coefficients exist for all $L_p$-functions. We let
$\d_p(\EuScript A(\mathbb{F}),\EuScript B(\mathbb{F}))$ denote the
smallest constant $c$ such that
\[ \Big\| \sum_{k\in\mathbb{F}} \sprod f{b_k} a_k \Big\|_p \leq c
\| f\|_p
\]
holds for all functions $f\in L_p[0,1]$. Observe that $\d_p(\EuScript
A(\mathbb{F}),\EuScript B(\mathbb{F}))$ is the norm of the operator
$S:L_p[0,1] \rightarrow L_p[0,1]$ given by $Sf = \sum_{k\in
  \mathbb{F}} \sprod f{b_k} a_k$.

We have the following facts about the quantities $\r_p$ and $\d_p$,
which are either obvious or proved in \cite{pietsch98:_orthon_banac}.

\begin{eqnarray}
  \label{eq:r_d}
  \r_p(\EuScript A(\mathbb{F}),\EuScript B(\mathbb{F}))
  & \leq &
  \d_p(\EuScript A(\mathbb{F}),\EuScript B(\mathbb{F})), \\
  \label{eq:dual}
  \d_p(\EuScript A(\mathbb{F}),\EuScript B(\mathbb{F}))
  & = &
  \d_{p'}(\EuScript B(\mathbb{F}),\EuScript A(\mathbb{F})), \\
  \label{eq:proj}
  \d_p(\EuScript A(\mathbb{F}),\EuScript B(\mathbb{F}))
  & \leq &
  \r_p(\EuScript A(\mathbb{F}),\EuScript B(\mathbb{F}))
  \d_p(\EuScript B(\mathbb{F}),\EuScript B(\mathbb{F})), \\
  \r_p(\EuScript A^\sigma(\mathbb{F}),\EuScript
  B^\sigma(\mathbb{F})) & = &
  \r_p(\EuScript A(\sigma(\mathbb{F})),\EuScript
  B(\sigma(\mathbb{F}))), \label{eq:1} \\
  \d_p(\EuScript A^\sigma(\mathbb{F}),\EuScript
  B^\sigma(\mathbb{F})) & = &
  \d_p(\EuScript A(\sigma(\mathbb{F})),\EuScript
  B(\sigma(\mathbb{F}))). \label{eq:2}
\end{eqnarray}
If $\mathbb{F}\subseteq\mathbb{G}$ then also
\begin{equation}
 \label{eq:monotone}
  \r_p(\EuScript A(\mathbb{F}),\EuScript B(\mathbb{F}))
  \leq
  \r_p(\EuScript A(\mathbb{G}),\EuScript B(\mathbb{G})).
\end{equation}
Moreover, if $\theta\in (0,1)$ is given by $1/p = (1-\theta)/p_0 +
\theta/p_1$, complex interpolation shows that
\begin{equation}
 \label{eq:interpol}
  \d_p(\EuScript A(\mathbb{F}),\EuScript B(\mathbb{F}))
  \leq
  \d_{p_0}(\EuScript A(\mathbb{F}),\EuScript B(\mathbb{F}))^{1-\theta}
  \d_{p_1}(\EuScript A(\mathbb{F}),\EuScript
  B(\mathbb{F}))^{\theta}.
\end{equation}

The next fact follows from the boundedness of the Riesz transform in
any $L_p[0,1]$ for $1<p<\infty$, see \cite[vol. I, p.
67]{zygmund59:_trigon}, and the boundedness of the canonical
projection from $L_p[0,1]$ onto the span of the first $2^m$ Walsh
functions, see \cite[p. 142]{schipp90:_walsh}.
\begin{lemma}
  For $1<p<\infty$, there is a constant $c_p$ such that
  \begin{equation}\label{eq:4}
    \d_p(\EuScript E_n,\EuScript E_n) \leq c_p
  \end{equation}
  for $n=1,2,\ldots.$ Moreover, we have for all $m=1,2,\dots$ that
  \begin{equation}\label{eq:3}
    \d_p(\EuScript W_{2^m},\EuScript W_{2^m}) = 1.
  \end{equation}
\end{lemma}
Combining  (\ref{eq:1}) and (\ref{eq:monotone}) gives the
following lemma.
\begin{lemma}
  If $\sigma[n]\subseteq[N]$ then
  \begin{equation}\label{eq:5}
    \r_p(\EuScript A_n^\sigma,\EuScript B_n) \leq \r_p(\EuScript A_N,
    \EuScript B_N^{\sigma^{-1}}).
  \end{equation}
\end{lemma}
\begin{proof}
  We have
  \begin{eqnarray*}
    \r_p(\EuScript A_n^\sigma,\EuScript B_n) & = &
    \r_p(\EuScript A(\sigma[n]),\EuScript
    B^{\sigma^{-1}}(\sigma[n])) \\
    & \leq & \r_p(\EuScript A([N]),\EuScript B^{\sigma^{-1}}([N]))
    =
    \r_p(\EuScript A_N, \EuScript B_N^{\sigma^{-1}}).
  \end{eqnarray*}
\end{proof}
Next we prove a first duality result.
\begin{lemma}
  For any orthonormal system $\EuScript A_n$ and $1<p<\infty$ we have
  \begin{eqnarray}
    \label{eq:8}
    \r_p(\EuScript E_n,\EuScript A_n)
    & \leq &
    c_p \r_{p'}(\EuScript A_n,\EuScript E_n),\\
    \r_p(\EuScript W_{2^m},\EuScript A_{2^m})
    & \leq &
    \r_{p'}(\EuScript A_{2^m},\EuScript W_{2^m}).\label{eq:9}
  \end{eqnarray}
\end{lemma}
\begin{proof}
  It follows successively from (\ref{eq:r_d}),
  (\ref{eq:dual}),~(\ref{eq:proj}) and (\ref{eq:4}) that
  \begin{eqnarray*}
    \r_p(\EuScript E_n,\EuScript A_n)
    & \leq & \d_p(\EuScript E_n,\EuScript A_n)
    = \d_{p'}(\EuScript A_n,\EuScript E_n)
    \leq \r_{p'}(\EuScript A_n,\EuScript E_n)
    \d_{p'}(\EuScript E_n,\EuScript E_n) \\
    & \leq &
    c_p \r_{p'}(\EuScript A_n,\EuScript E_n).
  \end{eqnarray*}
  The second inequality follows in the same way using (\ref{eq:3})
  instead of (\ref{eq:4}).
\end{proof}
We can now prove the complete duality result.
\begin{proposition}
  \label{prop:2}
  Given a permutation $\sigma$ and $n\in\{0,1,\dots\}$ there exists a
  number $N$ such that
  \begin{align}
    \r_p(\EuScript E_n,\EuScript W^\sigma_n) & \leq c_p
    \r_{p'}(\EuScript W_N^\sigma,\EuScript E_N),&
    \r_p(\EuScript E_n^{\sigma},\EuScript W_n) & \leq c_p
    \r_{p'}(\EuScript W_N,\EuScript E_N^{\sigma}), \label{eq:6}\\
    \r_p(\EuScript W_n,\EuScript E_n^{\sigma}) & \leq \phantom{c_p}
    \r_{p'}(\EuScript E_N^{\sigma},\EuScript W_N),&
    \r_p(\EuScript W_n^\sigma,\EuScript E_n) & \leq \phantom{c_p}
    \r_{p'}(\EuScript E_N,\EuScript W_N^\sigma).\label{eq:7}
  \end{align}
  for all $1<p<\infty$.
\end{proposition}
\begin{proof}
  The left hand inequalities are immediate consequences of
  (\ref{eq:8}) and (\ref{eq:9}) using $N=n$ and $N=2^m>n$ and
  (\ref{eq:monotone}) respectively.

  The right hand inequalities follow from (\ref{eq:5}), the
  corresponding left hand inequalities and (\ref{eq:5}) again.
\end{proof}
We can now summarize the duality results as follows.
\begin{proposition}
  \label{duality}
  Let $1<p<\infty$. Then the systems $(e_k)_{k\ge 0}$ and
  $(w_{\sigma(k)})_{k\ge 0}$ are equivalent in $L_p$ if and only if
  they are equivalent in $L_{p'}$.
\end{proposition}
\begin{proof}
  We only have to note that $(e_k)_{k\ge 0} \sim_p
  (w_{\sigma(k)})_{k\ge 0}$ if and only if the parameters
  $\r_p(\EuScript E_n,\EuScript W_n^\sigma)$ and $\r_p(\EuScript
  W_n^\sigma,\EuScript E_n)$ are uniformly bounded.
\end{proof}

%%%%%%%%%%%%%%%%%%%%%%%%%%%%%%%%%%%%%%%%%%%%%%%%%%%%%%%%%%%%%%%%%%%%%%%

\section{The key functions}

To show non-equivalence of the trigonometric and rearranged Walsh
system, norm estimates for the functions
\[ F^\sigma_n(s,t) = \sum_{k=0}^{n-1} e_k(s) w_{\sigma(k)}(t)
\]
play an essential r{\^o}le. This is due to the next observation.

\begin{proposition}
  \label{lower_bound}
  For each $p$ with $1< p < \infty$, there exists some constant
  $c_p>0$ such that
  \[ \r_p(\EuScript E_n,\EuScript W^\sigma_n)
  \ge c_p n^{1-1/p} \|F^\sigma_n\|_p^{-1} \ \ \mbox{for} \
  n=1,2,\ldots,
  \]
  where $\|F^\sigma_n\|_p$ is the norm of $F^\sigma_n$ in $L_p([0,1]^2)$.
\end{proposition}

\begin{proof}
  From the definition of $\r_p(\EuScript E_n,\EuScript W^\sigma_n)$ we
  find that
  \[ \int_0^1 \left| \sum_{k=0}^{n-1} \xi_k e_k(t) \right|^p dt \le
  \r_p(\EuScript E_n,\EuScript W^\sigma_n)^p
  \int_0^1 \left| \sum_{k=0}^{n-1} \xi_k w_{\sigma(k)}(t)
  \right|^p dt
  \]
  for all complex numbers $\xi_0,\ldots,\xi_{n-1}$. Using this for
  $\xi_k = e_k(s)$, integrating over $s\in [0,1]$ and taking $p$-th
  roots, we obtain that
  \[ \left( \int_0^1 \left| \sum_{k=0}^{n-1} e_k(t) \right|^p dt
  \right)^{1/p} \le \r_p(\EuScript E_n,\EuScript W^\sigma_n)
  \|F^\sigma_n\|_p.
  \]
  The left hand side is the $L_p$-norm of the Dirichlet kernel.  The
  well-known properties of this kernel imply that
  \[ c_1 n^{1-1/p} \le \left( \int_0^1 \left| \sum_{k=0}^{n-1} e_k(t)
    \right|^p dt \right)^{1/p} \le c_2 n^{1-1/p}
  \]
  where $c_1$ and $c_2$ depend only on $p$, see \cite[vol. I, p.
  67]{zygmund59:_trigon}. This completes the proof.
\end{proof}

Since by Parseval's equality $\|F^\sigma_n\|_2=\sqrt{n}$ and since
obviously $\|F^\sigma_n\|_\infty=n$, H{\"o}lder's inequality yields
the upper bound $\|F^\sigma_n\|_p\le n^{1-1/p}$ for any $2<p<\infty$.
If we can show for some $p>2$ that actually
\begin{equation}
  \label{crucialasymp}
  \liminf_{n\to\infty} n^{1/p-1} \|F^\sigma_n\|_p = 0
\end{equation}
then Proposition \ref{lower_bound} gives that $(e_k)_{k\ge 0}$ and
$(w_{\sigma(k)})_{k\ge 0}$ can not be equivalent in $L_p$. The duality
result Proposition \ref{duality} shows that this is also true in
$L_{p'}$.

We now derive that (\ref{crucialasymp}) for one $p_0$ in $(2,\infty)$
already implies (\ref{crucialasymp}) for all $p$ in $(2,\infty)$.
Indeed, if $p\in (2,p_0)$, defining $\theta\in (0,1)$ by $1/p=\theta/2
+ (1-\theta)/p_0$, H{\"o}lder's inequality together with
$\|F^\sigma_n\|_2=n^{1/2}$ yields
\[ n^{1/p-1} \|F^\sigma_n\|_p \le n^{1/p-1} \|F^\sigma_n\|_2^\theta
\|F^\sigma_n\|_{p_0}^{1-\theta} = \left( n^{1/p_0-1} \|F^\sigma_n\|_{p_0}
\right)^{1-\theta}.
\]
Similarly, if $p\in (p_0,\infty)$, defining $\theta\in (0,1)$ by $1/p
= \theta/p_0$, we obtain from $\|F^\sigma_n\|_\infty=n$ that
\[ n^{1/p-1} \|F^\sigma_n\|_p \le n^{1/p-1} \|F^\sigma_n\|_{p_0}^\theta
\|F^\sigma_n\|_\infty^{1-\theta} = \left( n^{1/p_0-1} \|F^\sigma_n\|_{p_0}
\right)^{\theta}.
\]

Altogether, we have proved the following theorem.
\begin{theorem}
  \label{thereisonlyonep}
  If there exists $p_0\in (2,\infty)$ such that
  \[ \liminf_{n\to\infty} n^{1/p_0-1} \|F^\sigma_n\|_{p_0} = 0
  \]
  then, for all $p\in (1,\infty)$ with $p\neq 2$, the systems $(e_k)$
  and $(w_{\sigma(k)})$ are not equivalent in $L_p$.
\end{theorem}

\begin{remark}
  In the cases $p=1$ and $p=\infty$, some additional care has to be
  taken. Using the well-known estimates
  \[ \d_1(\EuScript E_n,\EuScript E_n) =
  \d_\infty(\EuScript E_n,\EuScript E_n) \leq c (1+\log n),\]
  one obtains with the interpolation formula (\ref{eq:interpol})
  that if there exists $p_0\in (2,\infty)$ such that
  \[ \liminf_{n\to\infty} n^{1/p_0-1} (1+\log n)^{1-2/p_0}
  \|F^\sigma_n\|_{p_0} = 0
  \]
  then $(e_k)$ and $(w_{\sigma(k)})$ are not equivalent in $L_1$ and
  $L_\infty$.
\end{remark}

%%%%%%%%%%%%%%%%%%%%%%%%%%%%%%%%%%%%%%%%%%%%%%%%%%%%%%%%%%%%%%%%%%%%%%%

\section{The case $p=4$}

By the results of the previous section, we can now concentrate on
upper bounds for the $L_p$-norm of $F^\sigma_n$ for a convenient value
of $p\in (2,\infty)$. We use $p=4$ here. By Theorem
\ref{thereisonlyonep}, to show non-equivalence of $(e_k)$ and
$(w_{\sigma(k)})$ in $L_p$ for all $p\in(1,\infty)$ it is enough to
verify that
\[  \liminf_{n\to\infty} n^{-3/4} \|F^\sigma_n\|_4 = 0.
\]

We are going to formulate an equivalent combinatorial condition. To
this end, let us introduce some more notation. Given two numbers
$m,n\in \mathbb{N}_0$ with binary expansions $m=\sum_{i=0}^\infty m_i
2^i$ and $n=\sum_{i=0}^\infty n_i 2^i$, the {\em dyadic sum} is given
by $m\oplus n = \sum_{i=0}^\infty |m_i - n_i| 2^i$. The set
$\mathbb{N}_0$ with dyadic addition is isomorphic to the group of
Walsh functions which is expressed in the equation $w_{m\oplus n} =
w_m w_n$ for all $m,n \in \mathbb{N}_0$. Since the Walsh system claims
relationship with the powers of two, we will from now on concentrate
on the norms of $F_{2^n}$ instead of $F_n$ for all $n\in\mathbb{N}_0$.
Let
\[ A^\sigma_n = \{ (k,l,m) \in [2^n]^3 :k+l-m \in[2^n],
\sigma(k) \oplus \sigma(l) \oplus \sigma(m) = \sigma(k+l-m) \}.
\]

In the next lemma and throughout the paper, the notation $\#A$ means
the cardinality of a set $A$.

\begin{lemma}
  \label{computed4norm}
  $ \|F^\sigma_{2^n}\|_4^4 = \# A^\sigma_n$.
\end{lemma}

\begin{proof}
  It follows from
  \begin{eqnarray*}
    |F^\sigma_{2^n}(s,t)|^2 &=& F^\sigma_{2^n}(s,t)
    \overline{F^\sigma_{2^n}(s,t)} =
    \Big( \sum_{k=0}^{2^n-1} e_k(s) w_{\sigma(k)}(t) \Big)
    \Big( \sum_{l=0}^{2^n-1} e_{-l}(s) w_{\sigma(l)}(t) \Big) \\
    &=& \sum_{k,l=0}^{2^n-1} e_{k-l}(s) w_{\sigma(k)\oplus
    \sigma(l)}(t)
  \end{eqnarray*}
  that
  \[ |F^\sigma_{2^n}(s,t)|^4 = \sum_{k_1,l_1=0}^{2^n-1}
  \sum_{k_2,l_2=0}^{2^n-1}
  e_{k_1-l_1+k_2-l_2}(s) w_{\sigma(k_1)\oplus \sigma(l_1)
    \oplus \sigma(k_2)\oplus \sigma(l_2)}(t).
  \]
  Since
  \[ \int_0^1 \int_0^1 e_a(s) w_b(t) \,ds\,dt =
  \left\{ \begin{array}{ll}
      1 & \mbox{\ if}\  a=b=0, \\
      0 & \mbox{\ otherwise,}
    \end{array}\right.
  \]
  we obtain by integration that $ \|F^\sigma_{2^n}\|_4^4 = \#
  B_n^\sigma$ where
  \[
  B^\sigma_n =
  \Big\{ (k_1,l_1,k_2,l_2)\in [2^n]^4 :
  \begin{array}{l}
    k_1-l_1+k_2-l_2 = 0,\\
    \sigma(k_1)\oplus \sigma(l_1)
    \oplus \sigma(k_2)\oplus \sigma(l_2)=0
  \end{array}
  \Big\}.
  \]
  Obviously, $B^\sigma_n$ has the same cardinality as $A^\sigma_n$.
\end{proof}

\begin{corollary}
  \label{noneqcomb}
  If $\liminf_{n\to\infty} 8^{-n} \# A^\sigma_n = 0$ then, for all
  $p\in (1,\infty)$ with $p\neq 2$, the systems $(e_k)$ and
  $(w_{\sigma(k)})$ are not equivalent in $L_p$.
\end{corollary}

\begin{remark}
  Using the remark following Theorem \ref{thereisonlyonep} we also
  obtain that
  \[ \liminf_{n\to\infty} 8^{-n} n^2 \# A^\sigma_n = 0
  \]
  implies that $(e_k)$ and $(w_{\sigma(k)})$ are not equivalent in
  $L_1$ and $L_\infty$.
\end{remark}

%%%%%%%%%%%%%%%%%%%%%%%%%%%%%%%%%%%%%%%%%%%%%%%%%%%%%%%%%%%%%%%%%%%%%%%

\section{Application to concrete rearrangements}

In this section, we apply the results of the previous section to the
study of the equivalence problem for some specific rearrangements. In
particular, we treat the (besides the Walsh-Paley order) most
frequently used cases of the original Walsh system, the Walsh-Kaczmarz
system and the Walsh-Kronecker systems. For the properties and
alternative definitions of the above orderings, we refer the reader to
\cite{schipp90:_walsh}.

%%%%%%%%%%%%%%%%%%%%%%%%%%%%%%%%%%%%%%%%%%%%%%%%%%%%%%%%%%%%%%%%%%%%%%%

\subsection{Dyadically linear rearrangements}

The original Walsh system is a particular case of a linear
rearrangement of the Walsh-Paley system. A {\em dyadically linear
  rearrangement} is represented by a matrix
$T=(t_{i,j})_{i,j=0}^\infty$ with entries in $\{0,1\}$ such that the
$i$-th coefficient in the binary expansion of $\sigma(n)$ is given as
\[ \sigma(n)_i = \sum_{j=0}^\infty t_{i,j} n_j \mod 2.
\]
This is equivalent to the condition that $\sigma$ is linear with
respect to binary addition: $\sigma(m\oplus n) = \sigma(m) \oplus
\sigma(n)$. The {\em original Walsh system} is obtained using the
matrix $T$ with entries $t_{i,j}=1$ if and only if $j=i$ or $j=i+1$.

For linear rearrangements $\sigma$ the sets $A_n^\sigma$ behave
nicely.

\begin{proposition}\label{prop:1}
  If $\sigma$ is dyadically linear and $\pi$ is an arbitrary
  permutation, then
  \[ A_n^{\sigma\circ\pi} = A_n^\pi
  \]
  and consequently
  \[ \# A_n^{\sigma\circ\pi} = \# A_n^\pi.
  \]
\end{proposition}
\begin{proof}
  We simply observe that by linearity and injectivity of $\sigma$ we
  have
  \begin{eqnarray*}
    A_n^{\sigma\circ\pi}
    & = &
    \Big\{ (x,y,z)\in[2^n]^3 :
    \begin{array}{l}
      x+y-z \in[2^n],\\
      \sigma(\pi(x))\oplus\sigma(\pi(y))\oplus\sigma(\pi(z))
      = \sigma(\pi(x+y-z))
    \end{array}
    \Big\} \\
    & = &
    \Big\{ (x,y,z)\in[2^n]^3 :
    \begin{array}{l}
      x+y-z \in[2^n],\\
      \sigma(\pi(x)\oplus\pi(y)\oplus\pi(z))
      = \sigma(\pi(x+y-z))
    \end{array}
    \Big\} \\
    & = &
    \Big\{ (x,y,z)\in[2^n]^3 :
    \begin{array}{l}
      x+y-z \in[2^n],\\
      \pi(x)\oplus\pi(y)\oplus\pi(z)
      = \pi(x+y-z)
    \end{array}
    \Big\}
    = A_n^\pi.
  \end{eqnarray*}
\end{proof}

To use our general combinatorial condition for dyadically linear
rearrangements of the Walsh-Paley system, we need the following result
which may also have some interest in itself.

\begin{theorem}
  \label{combforlin}
  Let $\psi : \mathbb{N}_0 \rightarrow \mathbb{Z}$ be an arbitrary map.
  Then for all $n=0,1,\ldots$ we have
  \[ \# \{(x,y)\in[2^n]^2: \psi(x\oplus y)=x+y \}\leq3^n. \]
\end{theorem}
\begin{proof}
  For $u = 0,1,\ldots$, define
  \begin{eqnarray*}
    B_n(u)      &=& \{ (x,y) \in [2^n]^2 \ :\ x\oplus y = u\}
    \ \mbox{and}  \\
    C_n^\psi(u) &=& \{ (x,y) \in [2^n]^2 \ :\ x+y = \psi(u)\}.
  \end{eqnarray*}
  Then
  \[ \{(x,y)\in[2^n]^2\ : \  \psi (x\oplus y) = x+y \} = \bigcup_u
  B_n(u)\cap C_n^\psi(u).
  \]
  So all we have to show is
  \[ \sum_u \# \big(B_n(u) \cap C_n^\psi(u)\big) \le 3^n.\]

  We use induction over $n$. The statement for $n=0$ is trivial.  So
  assume we already know the statement for a certain value of $n$ and
  all functions $\psi$. Let us partition $B_{n+1}(u)$ into four
  disjoint subsets as follows
  \begin{eqnarray*}
    B_{00}(u) &=& B_{n+1}(u) \cap ([2^n] \times [2^n]), \\
    B_{01}(u) &=& B_{n+1}(u) \cap ([2^n] \times (2^n+[2^n])), \\
    B_{10}(u) &=& B_{n+1}(u) \cap ((2^n+[2^n]) \times [2^n]), \\
    B_{11}(u) &=& B_{n+1}(u) \cap ((2^n+[2^n]) \times (2^n+[2^n])).
  \end{eqnarray*}
  We are going to use the induction hypothesis to show that
  \begin{eqnarray}
    \label{first}
    \sum_u \# \big(B_{01}(u)                 \cap C_{n+1}^\psi(u)\big)
    &\le& 3^n, \\
    \label{second}
    \sum_u \# \big(B_{10}(u)                 \cap C_{n+1}^\psi(u)\big)
    &\le& 3^n, \\
    \label{third}
    \sum_u \# \big((B_{00}(u)\cup B_{11}(u)) \cap C_{n+1}^\psi(u)\big)
    &\le& 3^n.
  \end{eqnarray}
  This implies $\sum_u \# \big( B_{n+1}(u) \cap C_{n+1}^\psi(u)\big)
  \le 3^{n+1}$, completing the induction.

  To verify (\ref{first}), we observe that for $(x,y) \in [2^n] \times
  (2^n+[2^n])$ we have $y\oplus2^n=y-2^n$ and therefore
  \begin{eqnarray*}
    \lefteqn{
      (x,y) \in B_{01}(u) \cap C_{n+1}^\psi(u)
      \Longleftrightarrow x\oplus y = u
      \ \mbox{and}\ x+y=\psi(u)} \hspace{2cm} \\
    &\Longleftrightarrow& x\oplus y \oplus 2^n = u\oplus2^n
    \ \mbox{and}\ x+y-2^n=\psi(u)-2^n \\
    &\Longleftrightarrow& (x,y-2^n) \in B_n(u\oplus2^n) \cap
    C_{n+1}^{\tilde{\psi}}(u\oplus2^n),
  \end{eqnarray*}
  where we define $\tilde{\psi}(\tilde{u}) =
  \psi(\tilde{u}\oplus2^n)-2^n$. So
  \[ \#\big( B_{01}\cap
  C_{n+1}^\psi(u)\big) = \# \big( B_n(u\oplus2^n)\cap
  C_n^{\tilde{\psi}}(u\oplus2^n)\big)
  \]
  which yields by induction
  hypothesis that
  \[ \sum_u \# \big(B_{01}(u) \cap C_{n+1}^\psi(u)\big)
  \le\sum_u \# \big(B_n(u) \cap C_n^{\tilde{\psi}}(u)\big) \le 3^n.
  \]

  The inequality (\ref{second}) is symmetric to (\ref{first}).

  To prove (\ref{third}), we observe that
  \begin{eqnarray*}
    (x,y) \in B_{00}(u)  &\ &\mbox{implies}\ \ x+y <   2^{n+1} \
    \ \mbox{and} \\
    (x,y) \in B_{11}(u)  &\ &\mbox{implies}\ \ x+y \ge 2^{n+1}
  \end{eqnarray*}
  which gives that
  \begin{eqnarray*}
    \mbox{if}\quad \psi(u)\ge 2^{n+1} &\ &\mbox{then}\ \
    B_{00}(u) \cap C_{n+1}^\psi(u) = \emptyset \quad\mbox{and} \\
    \mbox{if}\quad \psi(u) < 2^{n+1}  &\ &\mbox{then}\ \
    B_{11}(u) \cap C_{n+1}^\psi(u) = \emptyset.
  \end{eqnarray*}
  So
  \begin{eqnarray*}
    \lefteqn{\sum_u \# \big((B_{00}(u)\cup B_{11}(u)) \cap
      C_{n+1}^\psi(u) \big)} \hspace{1cm}
    \\
    &=&
    \sum_{\psi(u)<2^{n+1}} \# \big(B_{00}(u) \cap C_{n+1}^\psi(u)\big)
    +
    \sum_{\psi(u)\ge 2^{n+1}} \# \big(B_{11}(u) \cap C_{n+1}^\psi(u)\big).
  \end{eqnarray*}
  Defining  $\tilde{\psi}$ by
  \[ \tilde{\psi}(u) = \left\{ \begin{array}{ll} \psi(u) & \
      \mbox{if}\
      \ \psi(u)<2^{n+1}, \\ \psi(u)-2^{n+1} & \ \mbox{if} \ \
      \psi(u)\ge 2^{n+1},
    \end{array} \right.
  \]
  we obtain for $u$ with $\psi(u) < 2^{n+1}$ that
  \[
  (x,y) \in B_{00}(u) \cap C_{n+1}^\psi(u)
  \Longleftrightarrow
  (x,y) \in B_n(u) \cap C_n^{\tilde{\psi}}(u)
  \]
  and for $u$ with $\psi(u)\ge 2^{n+1}$ that
  \[
  (x,y) \in B_{11}(u) \cap C_{n+1}^\psi(u)
  \Longleftrightarrow
  (x-2^n,y-2^n) \in B_n(u) \cap C_n^{\tilde{\psi}}(u).
  \]
  So
  \[ \sum_{\psi(u)<2^{n+1}} \# \big(B_{00}(u) \cap
  C_{n+1}^\psi(u)\big)
  =
  \sum_{\psi(u)<2^{n+1}} \# \big(B_n(u) \cap
  C_n^{\tilde{\psi}}(u) \big)
  \]
  and
  \[ \sum_{\psi(u)\geq 2^{n+1}} \# \big(B_{11}(u) \cap
  C_{n+1}^\psi(u)\big)
  =
  \sum_{\psi(u)\geq 2^{n+1}} \# \big(B_n(u) \cap
  C_n^{\tilde{\psi}}(u) \big)
  \]
  finally imply with the induction hypothesis that
  \[ \sum_u \# \big((B_{00}(u)\cup B_{11}(u)) \cap
  C_{n+1}^\psi(u)\big) \le
  \sum_u \# \big(B_n(u) \cap C_n^{\tilde{\psi}}(u)\big) \le 3^n.
  \]
\end{proof}

Denoting by $\iota$ the identity $\iota(x)=x$ for all
$x\in\mathbb{N}_0$ we can now prove the following result.
\begin{corollary}\label{cor:1}
  $\#A^\iota_n \le 6^n.$
\end{corollary}
\begin{proof}
  For each $z\in [2^n]$, we consider the set
  \[ A_n(z) = \{ (x,y)\in[2^n]^2 : x+y-z\in[2^n],
  x \oplus y \oplus z = x+y-z \}.
  \]
  Defining $\psi(u)= (u \oplus z) + z$, we obtain
  \[ A_n(z)\subseteq\{(x,y)\in[2^n]^2:\psi(x\oplus y)=x+y\},
  \]
  so from Theorem~\ref{combforlin} we infer that $\# A_n(z) \le 3^n$.
  Consequently
  \[ \#A^\iota_n \le \sum_{z\in [2^n]} \# A_n(z) \le 2^n
  3^n = 6^n.
  \]
\end{proof}

\begin{theorem}
  \label{noneqlinear}
  If $\sigma$ is dyadically linear then $ \#A^\sigma_n \le 6^n $ for
  $n=0,1,\ldots$. So the systems $(e_k)$ and $(w_{\sigma(k)})$ are not
  equivalent in $L_p$ for $p\neq 2$. In particular, the Walsh-Paley
  and the original Walsh system are not equivalent to the
  trigonometric system in $L_p$ for $p\in [1,\infty]$ with $p\neq 2$.
\end{theorem}

\begin{proof}
  The assertion follows immediately from Proposition~\ref{prop:1} and
  Corollaries~\ref{cor:1} and~\ref{noneqcomb} and the remark following
  Corollary~\ref{noneqcomb}.
\end{proof}

\begin{remark}
  The {\em Walsh-Kronecker systems} $W^{\sigma_n}_{2^n}$ are special
  rearrangements of the first $2^n$ Walsh functions different for each
  $n$ which are the basis for the fast Walsh-Fourier transform.  They
  can also be obtained from the Walsh matrices. Here $\sigma_n$ is a
  dyadically linear map on $[2^n]$ so that our results also apply to
  this case giving lower estimates for $\r_p(\EuScript
  E_{2^n},\EuScript W^{\sigma_n}_{2^n})$.
\end{remark}

%%%%%%%%%%%%%%%%%%%%%%%%%%%%%%%%%%%%%%%%%%%%%%%%%%%%%%%%%%%%%%%%%%%%%%%

\subsection{Piecewise linear rearrangements}

Unfortunately, one of the frequently used rearrangements of the Walsh
system, the {\em Walsh-Kaczmarz system}, is not a linear
rearrangement. It seems more natural in the equivalence problem than
the Walsh-Paley order since it arranges the Walsh functions in the
order of increasing number of sign changes. The corresponding
permutation $\sigma$ is given by $\sigma(0)=0$ and
\[ \sigma( 2^k + \sum_{i=0}^{k-1} x_i 2^i) = 2^k + \sum_{i=0}^{k-1}
x_{k-1-i} 2^i
\]
for $k=0,1,\ldots$ and $x_0,\ldots,x_{k-1}\in\{0,1\}$. It is possible
to estimate the cardinality of the set $A^\sigma_n$ from the previous
section for this rearrangement directly. Nevertheless, we prefer to
sketch an alternative approach which works for all piecewise linear
rearrangements.

A permutation $\sigma$ defines a {\em piecewise linear rearrangement}
if $\sigma(0)=0$ and
\[ \sigma( 2^k + m) = 2^k + \sigma_k(m)
\]
for $k=0,1,\ldots$, $0\le m\le 2^k-1$, and bijections $\sigma_k:[2^k]
\rightarrow [2^k]$ which are linear with respect to binary addition.
In particular, $\sigma$ leaves the blocks
$\{2^k,2^k+1,\ldots,2^{k+1}-1\}$ invariant. Obviously, the
Walsh-Kaczmarz order is a piecewise linear rearrangement.

Instead of using the functions $F^\sigma_n$, we now use the functions
\[ \tilde{F}^\sigma_n(s,t) = F^\sigma_{2n}(s,t) - F^\sigma_n(s,t) =
\sum_{k=n}^{2n-1}
e_k(s) w_{\sigma(k)}(t).
\]
As an analogue of Proposition \ref{lower_bound}, we obtain
\begin{proposition}
  \label{alt_lower_bound}
  For each $p$ with $1< p< \infty$, there exists some constant $c_p>0$
  such that
  \[ \r_p(\EuScript E_{2n},\EuScript W^\sigma_{2n})
  \ge c_p n^{1-1/p} \|\tilde{F}^\sigma_n\|_p^{-1} \ \ \mbox{for} \
  n=1,2,\ldots.
  \]
\end{proposition}

Similarly, we obtain analogues of Lemma \ref{computed4norm} and
Corollary \ref{noneqcomb} if we replace the set $A^\sigma_n$ by the
set
\[ \tilde{A}^\sigma_n =
\Big\{ (k,l,m) \in (2^n+[2^n])^3 :
\begin{array}{l}
  k+l-m \in 2^n+[2^n],\\
  \sigma(k) \oplus \sigma(l) \oplus \sigma(m) = \sigma(k+l-m)
\end{array}
\Big\}.
\]

\begin{lemma}
  \label{alt_computed4norm}
  $ \|\tilde{F}^\sigma_{2^n}\|_4^4 = \# \tilde{A}^\sigma_n$.
\end{lemma}

\begin{corollary}
  \label{alt_noneqcomb}
  If $ \liminf_{n\to\infty} 8^{-n} \# \tilde{A}^\sigma_n = 0$ then,
  for all $p\in (1,\infty)$ with $p\neq 2$, the systems $(e_k)$ and
  $(w_{\sigma(k)})$ are not equivalent in $L_p$.
\end{corollary}

\begin{remark}
  Using the remark following Theorem \ref{thereisonlyonep} we also
  obtain that
  \[ \liminf_{n\to\infty} 8^{-n} n^2 \# \tilde{A}^\sigma_n = 0
  \]
  implies that $(e_k)$ and $(w_{\sigma(k)})$ are not equivalent in
  $L_1$ and $L_\infty$.
\end{remark}

We are now in a position to treat the case of piecewise linear
rearrangements.

\begin{theorem}
  \label{noneqpiecelinear}
  If $\sigma$ is a piecewise linear rearrangement then $
  \#\tilde{A}^\sigma_n \le 6^n$ for $n=0,1,\ldots$. So the systems
  $(e_k)$ and $(w_{\sigma(k)})$ are not equivalent in $L_p$ for $p\neq
  2$. In particular, the Walsh-Kaczmarz system is not equivalent to
  the trigonometric system in $L_p$ for $p\in [1,\infty]$ with $p\neq
  2$.
\end{theorem}

\begin{proof}
  For each $z\in 2^n+[2^n]$, we consider the set
  \[ \tilde{A}^\sigma_n(z) =
  \Big\{ (x,y)\in (2^n+[2^n])^2 :
  \begin{array}{l}
    x+y - z \in 2^n+[2^n],\\
    \sigma(x) \oplus \sigma(y) \oplus \sigma(z) = \sigma(x+y-z) \}
  \end{array}
  \Big\}.
  \]
  Let $\sigma_k:[2^k] \rightarrow [2^k]$ denote the linear maps from
  the definition of piecewise linearity.  Then we obtain for any
  $(x,y)\in \tilde{A}^\sigma_n(z)$ that $x+y-z\in 2^n+[2^n]$ and
  \[ \sigma_n(\tilde{x}\oplus \tilde{y}\oplus \tilde{z})=
  \sigma_n(\tilde{x}+\tilde{y}-\tilde{z}),
  \]
  where $\tilde{x}=x-2^n, \tilde{y}-2^n,\tilde{z}=z-2^n$. Since
  $\sigma_n$ is a permutation this implies $\tilde{x}\oplus
  \tilde{y}\oplus \tilde{z} = \tilde{x}+\tilde{y}-\tilde{z}$.  So
  \[ \# \tilde{A}^\sigma_n(z) \le \# \{ (\tilde{x},\tilde{y}) \in [2^n]^2 :
  \tilde{x}\oplus \tilde{y}\oplus \tilde{z} =
  \tilde{x}+\tilde{y}-\tilde{z} \}.
  \]
  This can be estimated by $3^n$ as in the proof of Corollary
  \ref{cor:1} and gives $\#\tilde{A}^\sigma_n \le \sum_{z\in
    2^n+[2^n]} \# \tilde{A}^\sigma_n(z) \le 6^n$. The claim now
  follows from Corollary~\ref{alt_noneqcomb} if $p\in(1,\infty)$ and
  the remark after that corollary if $p=1,\infty$.
\end{proof}

%%%%%%%%%%%%%%%%%%%%%%%%%%%%%%%%%%%%%%%%%%%%%%%%%%%%%%%%%%%%%%%%%%%%%%%

\subsection{Small perturbations}
\label{sec:perturb}

Besides (piecewise) linear rearrangements, we can treat a further
class of rearrangements, namely small perturbations of rearrangements
of the Walsh system, that are known to be non-equivalent.

To this end, for $v\in\mathbb{Z}$, we also consider the sets
\[ A_n^\sigma(v) := \Big\{ (x,y,z) \in [2^n]^3 :
\begin{array}{l}
x+y-z\in[2^n],\\
\sigma(x) \oplus \sigma(y) \oplus \sigma(z) \oplus \sigma(v) =
\sigma(x + y - z)
\end{array}
\Big\},
\]
So instead of asking for $\sigma(x) \oplus \sigma(y) \oplus \sigma(z)
\oplus \sigma(x + y - z) =0$ we require that the left hand side of
this equality equals a fixed number $\sigma(v)$. Note that
$A_n^\sigma(0) = A_n^\sigma$.

As in the proof of Theorem \ref{noneqlinear}, we can control the size
of $A_n^\sigma(v)$ for dyadically linear rearrangements $\sigma$ and
all $v$.
\begin{proposition}
  If $\sigma$ is dyadically linear then $ \#A_n^\sigma(v) \le 6^n $
  for $n=0,1,\ldots$.
\end{proposition}
\begin{proof}
  By linearity and injectivity of $\sigma$ we need only consider the
  case $\sigma=\iota$. Using $\psi(u)= (u\oplus z\oplus v)+z$ we see
  the result as in the proof of Corollary~\ref{cor:1}.
\end{proof}

Given two permutations $\pi$ and $\sigma$ let
\[ f(u) := \pi(u)\oplus\sigma(u)
\]
and put
\[ f_n^* := \max_{u\in[2^n]} f(u).
\]
The function $f$ measures in some sense, how much $\pi$ deviates from
$\sigma$. In particular $|\pi(u)-\sigma(u)|\leq f(u)$. We say that
$\pi$ \emph{dyadically differs from} $\sigma$ by $f$.

\begin{proposition}
  \label{thisone}
  We have
  \[ A_n^\sigma \subseteq \bigcup_{\pi(v)\leq4f_n^*} A_n^\pi (v).
  \]
  In particular
  \[ \#A_n^\sigma \leq (4f_n^*+1)\,\#A_n^\pi.
  \]
\end{proposition}
\begin{proof}
  Note that for all $x,y,z$ we have
  \begin{multline}\label{eq:10}
    \sigma(x)\oplus\sigma(y)\oplus\sigma(z)\oplus\sigma(x+y-z)
    \oplus\pi(x)\oplus\pi(y)\oplus\pi(z)\oplus\pi(x+y-z) \\
    = f(x)\oplus f(y)\oplus f(z)\oplus f(x+y-z).
  \end{multline}
  Also, for any $x,y\geq0$ the dyadic addition satisfies
  \[ x\oplus y \leq x+y.
  \]
  Therefore if $(x,y,z)\in A_n^\sigma$ then $x,y,z,x+y-z\in[2^n]$ and
  hence
  \[ f(x)\oplus f(y)\oplus f(z)\oplus f(x+y-z) \leq 4f_n^*.
  \]
  Defining $v$ by
  \[
  \pi(v)=\pi(x)\oplus\pi(y)\oplus\pi(z)\oplus\pi(x+y-z),
  \]
  we obtain $ (x,y,z)\in A_n^\pi(v).  $ It now follows from
  \[ \sigma(x)\oplus\sigma(y)\oplus\sigma(z)\oplus\sigma(x+y-z) = 0
  \]
  and~(\ref{eq:10}) that $\pi(v) \leq 4f_n^*.$ This completes the
  proof.
\end{proof}

This proposition immediately implies

\begin{theorem}
  If $\pi$ dyadically differs from $\sigma$ by $f$, and $f$ and
  $\#A_n^\pi$ satisfy
  \[ \liminf_{n\to\infty} 8^{-n}\#A_n^\pi f_n^* = 0
  \]
  then, for all $p\in(1,\infty)$ with $p\not=2$, the systems $(e_k)$
  and $(w_{\sigma(k)})$ are not equivalent in $L_p$.  This is in
  particular the case, if $\pi$ is dyadically linear and $f$ satisfies
  \[ f_n^* = o\Big(\frac{4^n}{3^n}\Big).
  \]
\end{theorem}

\begin{remark}
  In the cases $p=1,\infty$, we again have to adjust the condition to
  \[ \liminf_{n\to\infty} 8^{-n}n^2\#A_n^\pi f_n^* = 0.
  \]
\end{remark}

We now develop a dual version of the last results. We will mostly
leave the proofs to the reader, since they are completely analogous to
the previous ones. For $v\in\mathbb{Z}$, we define the sets
\[ \hat A_n^\sigma(v) := \{ (x,y,z) \in [2^n]^3 :
\sigma(x) + \sigma(y) - \sigma(z) + v = \sigma(x \oplus y \oplus
z)\}
\]
and we let $\hat A_n^\sigma=\hat A_n^\sigma(0)$. As before, we can
show that for $p\in(1,\infty)$ with $p\not=2$ the systems
$(e_{\sigma(k)})$ and $(w_k)$ are not equivalent in $L_p$ if
\[ \liminf_{n\to\infty}8^{-n} \#\hat A_n^\sigma = 0.
\]

Again, we can control the size of $\hat A_n^\sigma(v)$ for dyadically
linear rearrangements $\sigma$ and all $v$.
\begin{proposition}
  If $\sigma$ is dyadically linear then $ \#\hat A_n^\sigma(v) \le 6^n
  $ for $n=0,1,\ldots$.
\end{proposition}

Given two permutations $\pi$ and $\sigma$ let
\[ \hat f(u) := |\pi(u)-\sigma(u)|
\]
and put
\[ \hat f_n^* := \max_{u\in[2^n]} \hat f(u).
\]
The function $\hat f$ measures how much $\pi$ deviates from $\sigma$.
We say that $\pi$ \emph{differs from} $\sigma$ by $\hat f$.

\begin{proposition}\label{prop:dual_variation}
  We have
  \[ \hat A_n^\sigma \subseteq \bigcup_{|v|\leq4\hat f_n^*}
  \hat A_n^\pi (v).
  \]
  In particular
  \[ \#\hat A_n^\sigma \leq (8\hat f_n^*+1) \, \#\hat A_n^\pi.
  \]
\end{proposition}
\begin{proof}
  Note that for all $x,y,z$ we have
  \begin{multline}\label{eq:11}
    |\sigma(x\oplus y\oplus z) - \sigma(x) - \sigma(y) + \sigma(z) -
    \pi(x\oplus y\oplus z) + \pi(x) + \pi(y) - \pi(z)| \\
    \leq \hat f(x\oplus y\oplus z) + \hat f(x) + \hat f(y) + \hat
    f(z).
  \end{multline}
  Therefore if $(x,y,z)\in A_n^\sigma$ then $x,y,z,x\oplus y\oplus
  z\in[2^n]$ and hence
  \[ \hat f(x\oplus y\oplus z) + \hat f(x) + \hat f(y) + \hat
  f(z) \leq 4\hat f_n^*.
  \]
  Defining $v$ by
  \[
  v=\pi(x\oplus y\oplus z) - \pi(x) - \pi(y) + \pi(z),
  \]
  we obtain $ (x,y,z)\in \hat A_n^\pi(v).  $ It now follows from
  \[ \sigma(x\oplus y\oplus z)- \sigma(x) -\sigma(y) +\sigma(z) = 0
  \]
  and~(\ref{eq:11}) that $ |v| \leq 4\hat f_n^*.  $ This completes the
  proof.
\end{proof}

Again we immediately obtain

\begin{theorem}
  If $\pi$ differs from $\sigma$ by $\hat f$, and $\hat f$ and $\#\hat
  A_n^\pi$ satisfy
  \[ \liminf_{n\to\infty} 8^{-n}\#\hat A_n^\pi \hat f_n^* = 0
  \]
  then, for all $p\in(1,\infty)$ with $p\not=2$, the systems
  $(e_{\sigma(k)})$ and $(w_k)$ are not equivalent in $L_p$.  This is
  in particular the case, if $\pi$ is dyadically linear and $\hat f$
  satisfies
  \[ \hat f_n^* = o\Big(\frac{4^n}{3^n}\Big).
  \]
\end{theorem}

\begin{remark}
  In the cases $p=1,\infty$, we again have to adjust the condition to
  \[ \liminf_{n\to\infty} 8^{-n}n^2\#\hat A_n^\pi \hat f_n^* = 0.
  \]
\end{remark}

To illustrate the power of this perturbation method, we add another
example.
\begin{example}
  Let $\mathbb F$ be a subset of $\mathbb{N}_0$ such that
  \[
  \liminf_{n\to \infty} \frac{\# ( \mathbb{F}\cap [n])}{n} <
  2- \log_2 3 = 0.415037\ldots .
  \]
  Let $\sigma$ be such that for $x=\sum_{i=0}^{\infty} x_i2^i$ we have
  \[ \sigma(x)
  = \sum_{i\in\mathbb F} \tilde x_i 2^i \oplus \sum_{i\not\in\mathbb
    F} x_i2^i,
  \]
  where $\tilde x_i\in\{0,1\}$ are such that $\sigma$ is a permutation
  and otherwise arbitrary. In other words, $\sigma$ acts arbitrarily
  on the binary coefficients in $\mathbb F$ and as the identity on the
  remaining binary coefficients.

  Then the systems $(e_k)$ and $(w_{\sigma(k)})$ are not equivalent in
  $L_p$ with $p\not=2$.
\end{example}
\begin{proof}
  Fix $n\in \mathbb{N}_0$.  Let $m= \#(\mathbb{F}\cap [n])$ and write
  $\mathbb F \cap [n] =\{k_0,\dots,k_{m-1}\}$ and $[n]\setminus\mathbb
  F=\{k_m,\dots,k_{n-1}\}$. Define a permutation $\pi_n$ by
  \[ \pi_n\Big(\sum_{i=0}^{n-1} x_i2^i\Big)
  = \sum_{i=0}^{n-1} x_{k_i} 2^i.
  \]
  Then $\pi_n$ is a dyadically linear permutation on $[2^n]$ so by
  Proposition~\ref{prop:1} we have
  $\#A^{\pi_n}_n=\#A^{\pi_n\circ\sigma}_n$. Moreover
  \[ \pi_n\circ\sigma(u) \oplus \pi_n(u) =  \sum_{i=0}^{m-1} \tilde
  u_{k_i}2^i \oplus \sum_{i=0}^{m-1} u_{k_i}2^i  < 2^m.
  \]
  This implies by Proposition~\ref{thisone} that
  \[ \#A_n^{\pi_n\circ\sigma} \leq 4\cdot 2^m \# A_n^{\pi_n}
  \]
  or by the linearity of $\pi_n$ and Corollary~\ref{cor:1}
  \[ \#A_n^\sigma \leq 4\cdot 2^m 6^n.
  \]
  The growth condition on $\# (\mathbb{F} \cap [n])$ ensures that
  \[ \liminf_{n\to\infty}8^{-n}\#A_n^\sigma = 0.
  \]
  The claim now follows from Corollary~\ref{noneqcomb}.
\end{proof}

\noindent{\bf Final remarks:} {\bf 1.} The estimates for the
non-equivalence quantities obtained by our methods have power type
behavior. Nevertheless, since they do not give optimal exponents
except possibly in the case $p=4$, we did not state those estimates
explicitly. In the case $2 < p\le 4$, our estimates for the
Walsh-Paley system are the same as the lower bounds obtained in
\cite{wojtaszczyk00:_non_walsh}. In the cases $p>4$ and $1\le p<2$,
the estimates for the special case of the Walsh-Paley order in
\cite{wojtaszczyk00:_non_walsh} are better than ours. It would be
interesting to find the optimal estimates at least in the cases of the
usual orderings.

{\bf 2.} Although we were not able to give general estimates for the
cardinalities of the sets $A_n^\sigma$, we conjecture that the
identical permutation already gives the maximal possible cardinality.
A similar and from the combinatorial point of view very natural
question is to find good upper bounds for the cardinalities of the
sets
\[ B^\sigma_n = \{ (k,l) \in [2^n]^2 :k+l \in[2^n],
\sigma(k) \oplus \sigma(l) = \sigma(k+l) \}.
\]
For linear and piecewise linear rearrangements one can obtain that $\#
B^\sigma_n\le 3^n$ and for the identity $\# B^\iota_n= 3^n$. Again we
conjecture that $\# B^\sigma_n\le 3^n$ holds for any permutation
$\sigma$. Basically, this is a question about how big the set of pairs
$(k,l)$ can be for which $\sigma$ behaves like a homomorphism between
the integers and the Cantor group. We checked this claim for $n\le 4$
and for all permutations $\sigma$ of $[2^n]$ with a computer. The
running time for the case $n=4$ on a PC was about four days.

%\bibliographystyle{abbrv}
%\bibliography{equi}

\noindent  Aicke Hinrichs, Mathematisches Institut, FSU Jena,
D-07743 {\sc Jena}, Germany, e-mail: {\tt nah@rz.uni-jena.de}

\medskip

\noindent J{\"o}rg Wenzel, Department of Mathematics and Applied
Mathematics, University of Pretoria, {\sc Pretoria} 0002, South
Africa,\\ e-mail: {\tt wenzel@minet.uni-jena.de}

\end{document}